\documentclass[12]{article}
\usepackage{bm}
\usepackage{fleqn}
\usepackage{graphicx}
\begin{document}
\Large
\begin{center}
\bf{The Complement of Binary Klein Quadric as\\ a Combinatorial Grassmannian}
\end{center}
\vspace*{-.5cm}
\begin{center}
Metod Saniga
\end{center}
\vspace*{-.5cm} \normalsize
\begin{center}
Institute for Discrete Mathematics and Geometry,
Vienna University of Technology\\ Wiedner Hauptstra\ss e 8--10,
A-1040 Vienna, Austria\\
(metod.saniga@tuwien.ac.at)

\vspace*{-.0cm}

and

\vspace*{-.0cm}

Astronomical Institute, Slovak Academy of Sciences\\
SK-05960 Tatransk\' a Lomnica, Slovak Republic\\
(msaniga@astro.sk)  
\end{center}

\vspace*{-.4cm} \noindent \hrulefill

\vspace*{-.1cm} \noindent {\bf Abstract}

\noindent 
Given a hyperbolic quadric of PG(5,\,2), there are 28 points off this quadric and 56 lines skew to it. It is shown that the $(28_6, 56_3)$-configuration formed by these points and lines is isomorphic to the combinatorial Grassmannian of type $G_2(8)$. It is also pointed out that a set of seven points of $G_2(8)$ whose labels share a mark corresponds to a Conwell heptad of PG(5,\,2). Gradual removal of Conwell heptads from
the $(28_6, 56_3)$-configuration yields a nested sequence of binomial configurations identical with part of that found to be associated with Cayley-Dickson algebras (arXiv:1405.6888).

\vspace*{.2cm}

\vspace*{-.2cm} \noindent \hrulefill

\vspace*{.1cm}
%\large
\noindent
Let $\mathcal{Q}^{+}(5,2)$ be a hyperbolic quadric in a five-dimensional projective space PG$(5,2)$. As it is well known (see, e.\,g., \cite{hir1,hir2}), there are 28 points lying off this quadric as well as 56 lines skew (or, external) to it. If the equation of the quadric is taken in a canonical form $x_1x_2 + x_3x_4 + x_5x_6 = 0$,
then the 28 off-quadric points are

\begin{center}
\begin{tabular}{||c|cccccc||}
\hline \hline
 No. & ~$x_1$ & $x_2$ & $x_3$ & $x_4$ & $x_5$ & $x_6$  \\
\hline
~1~  &    1   &   1   &   1   &   0   &   0   &   0    \\
~2~  &    1   &   1   &   0   &   0   &   1   &   0    \\
~3~  &    1   &   1   &   0   &   0   &   0   &   1    \\
~4~  &    1   &   1   &   0   &   1   &   0   &   0    \\
~5~  &    1   &   1   &   1   &   0   &   1   &   0    \\
~6~  &    1   &   1   &   1   &   0   &   0   &   1    \\
~7~  &    1   &   1   &   0   &   1   &   1   &   0    \\
~8~  &    1   &   1   &   0   &   1   &   0   &   1    \\
~9~  &    0   &   0   &   1   &   1   &   0   &   0    \\
~10~ &    0   &   0   &   1   &   1   &   1   &   0    \\
~11~ &    0   &   0   &   1   &   1   &   0   &   1    \\
~12~ &    0   &   1   &   1   &   1   &   0   &   0    \\
~13~ &    1   &   0   &   1   &   1   &   1   &   0    \\
~14~ &    1   &   0   &   1   &   1   &   0   &   1    \\
~15~ &    0   &   0   &   0   &   0   &   1   &   1    \\
~16~ &    1   &   0   &   0   &   0   &   1   &   1    \\
~17~ &    0   &   0   &   1   &   0   &   1   &   1    \\
~18~ &    0   &   0   &   0   &   1   &   1   &   1    \\
~19~ &    0   &   1   &   1   &   0   &   1   &   1    \\
~20~ &    0   &   1   &   0   &   1   &   1   &   1    \\
~21~ &    1   &   1   &   1   &   1   &   1   &   1    \\
\hline
~22~ &    1   &   1   &   0   &   0   &   0   &   0    \\
~23~ &    1   &   0   &   1   &   1   &   0   &   0    \\
~24~ &    0   &   1   &   1   &   1   &   1   &   0    \\
~25~ &    0   &   1   &   1   &   1   &   0   &   1    \\
~26~ &    0   &   1   &   0   &   0   &   1   &   1    \\
~27~ &    1   &   0   &   1   &   0   &   1   &   1    \\
~28~ &    1   &   0   &   0   &   1   &   1   &   1    \\
 \hline \hline
\end{tabular}
\end{center}
and the 56 external lines read

\begin{center}
\resizebox{\columnwidth}{!}{% 
{
\begin{tabular}{||c|ccccccccccccccccccccc|ccccccc||}
\hline \hline
 No.   & 1 & 2 & 3 & 4 & 5 & 6 & 7 & 8 & 9 & 10 & 11 & 12 & 13 & 14 & 15 & 16 & 17 & 18 & 19 & 20 & 21 & 22 & 23 & 24 & 25 & 26 & 27 & 28\\
\hline
~1~   & + &   &   & + &   &   &   &   & + &    &    &    &    &    &    &    &    &    &    &    &    &    &    &    &    &    &    &   \\
~2~   & + &   &   &   &   &   & + &   &   &  + &    &    &    &    &    &    &    &    &    &    &    &    &    &    &    &    &    &   \\ 
~3~   & + &   &   &   &   &   &   & + &   &    &  + &    &    &    &    &    &    &    &    &    &    &    &    &    &    &    &    &   \\
~4~   & + &   &   &   &   &   &   &   &   &    &    &    &    &    &    & +  &    &    &  + &    &    &    &    &    &    &    &    &   \\
~5~   & + &   &   &   &   &   &   &   &   &    &    &    &    &    &    &    &    & +  &    &    &  + &    &    &    &    &    &    &   \\
~6~   &   & + & + &   &   &   &   &   &   &    &    &    &    &    & +  &    &    &    &    &    &    &    &    &    &    &    &    &   \\
~7~   &   & + &   &   &   & + &   &   &   &    &    &    &    &    &    &    & +  &    &    &    &    &    &    &    &    &    &    &   \\
~8~   &   & + &   &   &   &   &   & + &   &    &    &    &    &    &    &    &    & +  &    &    &    &    &    &    &    &    &    &   \\
~9~   &   & + &   &   &   &   &   &   &   &    &  + &    &    &    &    &    &    &    &    &    &  + &    &    &    &    &    &    &   \\
~10~  &   & + &   &   &   &   &   &   &   &    &    & +  & +  &    &    &    &    &    &    &    &    &    &    &    &    &    &    &   \\
~11~  &   &   & + &   & + &   &   &   &   &    &    &    &    &    &    &    & +  &    &    &    &    &    &    &    &    &    &    &   \\
~12~  &   &   & + &   &   &   & + &   &   &    &    &    &    &    &    &    &    & +  &    &    &    &    &    &    &    &    &    &   \\
~13~  &   &   & + &   &   &   &   &   &   & +  &    &    &    &    &    &    &    &    &    &    &  + &    &    &    &    &    &    &   \\
~14~  &   &   & + &   &   &   &   &   &   &    &    & +  &    & +  &    &    &    &    &    &    &    &    &    &    &    &    &    &   \\
~15~  &   &   &   & + & + &   &   &   &   & +  &    &    &    &    &    &    &    &    &    &    &    &    &    &    &    &    &    &   \\
~16~  &   &   &   & + &   & + &   &   &   &    &  + &    &    &    &    &    &    &    &    &    &    &    &    &    &    &    &    &   \\
~17~  &   &   &   & + &   &   &   &   &   &    &    &    &    &    &    &  + &    &    &    &  + &    &    &    &    &    &    &    &   \\
~18~  &   &   &   & + &   &   &   &   &   &    &    &    &    &    &    &    &  + &    &    &    &  + &    &    &    &    &    &    &   \\
~19~  &   &   &   &   & + & + &   &   &   &    &    &    &    &    & +  &    &    &    &    &    &    &    &    &    &    &    &    &   \\
~20~  &   &   &   &   & + &   & + &   & + &    &    &    &    &    &    &    &    &    &    &    &    &    &    &    &    &    &    &   \\
~21~  &   &   &   &   & + &   &   &   &   &    &    &    &    & +  &    &    &    &    &    &  + &    &    &    &    &    &    &    &   \\
~22~  &   &   &   &   &   & + &   & + & + &    &    &    &    &    &    &    &    &    &    &    &    &    &    &    &    &    &    &   \\
~23~  &   &   &   &   &   & + &   &   &   &    &    &    & +  &    &    &    &    &    &    &  + &    &    &    &    &    &    &    &   \\
~24~  &   &   &   &   &   &   & + & + &   &    &    &    &    &    & +  &    &    &    &    &    &    &    &    &    &    &    &    &   \\
~25~  &   &   &   &   &   &   & + &   &   &    &    &    &    &  + &    &    &    &    &  + &    &    &    &    &    &    &    &    &   \\
~26~  &   &   &   &   &   &   &   & + &   &    &    &    &  + &    &    &    &    &    &  + &    &    &    &    &    &    &    &    &   \\
~27~  &   &   &   &   &   &   &   &   & + &    &    &    &    &    &    &    &  + &  + &    &    &    &    &    &    &    &    &    &   \\
~28~  &   &   &   &   &   &   &   &   & + &    &    &    &    &    &    &    &    &    &  + &  + &    &    &    &    &    &    &    &   \\
~29~  &   &   &   &   &   &   &   &   &   & +  & +  &    &    &    & +  &    &    &    &    &    &    &    &    &    &    &    &    &   \\
~30~  &   &   &   &   &   &   &   &   &   & +  &    &    &    &  + &    &  + &    &    &    &    &    &    &    &    &    &    &    &   \\
~31~  &   &   &   &   &   &   &   &   &   &    & +  &    & +  &    &    &  + &    &    &    &    &    &    &    &    &    &    &    &   \\
~32~  &   &   &   &   &   &   &   &   &   &    &    &  + &    &    &    &  + &    &    &    &    &  + &    &    &    &    &    &    &   \\
~33~  &   &   &   &   &   &   &   &   &   &    &    &  + &    &    &    &    &  + &    &    &  + &    &    &    &    &    &    &    &   \\
~34~  &   &   &   &   &   &   &   &   &   &    &    &  + &    &    &    &    &    &  + & +  &    &    &    &    &    &    &    &    &   \\
~35~  &   &   &   &   &   &   &   &   &   &    &    &    &  + &  + &  + &    &    &    &    &    &    &    &    &    &    &    &    &   \\
\hline
~36~  & + &   &   &   &   &   &   &   &   &    &    &    &    &    &    &    &    &    &    &    &    &    &    &    &    & +  & +  &   \\
~37~  &   & + &   &   &   &   &   &   &   &    &    &    &    &    &    &    &    &    &    &    &    &    &  + &  + &    &    &    &   \\
~38~  &   &   & + &   &   &   &   &   &   &    &    &    &    &    &    &    &    &    &    &    &    &    &  + &    &  + &    &    &   \\
~39~  &   &   &   & + &   &   &   &   &   &    &    &    &    &    &    &    &    &    &    &    &    &    &    &    &    & +  &    & + \\
~40~  &   &   &   &   & + &   &   &   &   &    &    &    &    &    &    &    &    &    &    &    &    &    &    &    &  + &    &    & + \\
~41~  &   &   &   &   &   & + &   &   &   &    &    &    &    &    &    &    &    &    &    &    &    &    &    &  + &    &    &    & + \\
~42~  &   &   &   &   &   &   & + &   &   &    &    &    &    &    &    &    &    &    &    &    &    &    &    &    &  + &    &  + &   \\
~43~  &   &   &   &   &   &   &   & + &   &    &    &    &    &    &    &    &    &    &    &    &    &    &    &  + &    &    &  + &   \\
~44~  &   &   &   &   &   &   &   &   & + &    &    &    &    &    &    &    &    &    &    &    &    &    &    &    &    &    &  + & + \\
~45~  &   &   &   &   &   &   &   &   &   &  + &    &    &    &    &    &    &    &    &    &    &    &    &    &    &  + &  + &    &   \\
~46~  &   &   &   &   &   &   &   &   &   &    &  + &    &    &    &    &    &    &    &    &    &    &    &    &  + &    &  + &    &   \\
~47~  &   &   &   &   &   &   &   &   &   &    &    &  + &    &    &    &    &    &    &    &    &    &  + &  + &    &    &    &    &   \\
~48~  &   &   &   &   &   &   &   &   &   &    &    &    &  + &    &    &    &    &    &    &    &    &  + &    &  + &    &    &    &   \\
~49~  &   &   &   &   &   &   &   &   &   &    &    &    &    &  + &    &    &    &    &    &    &    &  + &    &    &  + &    &    &   \\
~50~  &   &   &   &   &   &   &   &   &   &    &    &    &    &    &  + &    &    &    &    &    &    &    &    &  + &  + &    &    &   \\
~51~  &   &   &   &   &   &   &   &   &   &    &    &    &    &    &    &  + &    &    &    &    &    &  + &    &    &    &  + &    &   \\
~52~  &   &   &   &   &   &   &   &   &   &    &    &    &    &    &    &    &  + &    &    &    &    &    &  + &    &    &    &    & + \\
~53~  &   &   &   &   &   &   &   &   &   &    &    &    &    &    &    &    &    &  + &    &    &    &    &  + &    &    &    & +  &   \\
~54~  &   &   &   &   &   &   &   &   &   &    &    &    &    &    &    &    &    &    &  + &    &    &  + &    &    &    &    & +  &   \\
~55~  &   &   &   &   &   &   &   &   &   &    &    &    &    &    &    &    &    &    &    &  + &    &  + &    &    &    &    &    & + \\
~56~  &   &   &   &   &   &   &   &   &   &    &    &    &    &    &    &    &    &    &    &    &  + &    &  + &    &    &  + &    &   \\
 \hline \hline
\end{tabular}}%
} 
\end{center}

\bigskip
\noindent
In the latter table, the `+' symbol indicates which point lies on a given line; for example, line 1 consists of points 1, 4 and 9. As it is obvious from this table, each line has three points and through each point there are six lines; hence, these points and lines form a
$(28_6, 56_3)$-configuration.

Next, a combinatorial Grassmannian $G_k(X)$ (see, e.\,g., \cite{praz,owpr}), where $k$ is a positive
integer and $X$ is a finite set, $|X| =
N$, is a point-line incidence structure whose points are $k$-element subsets of $X$ and whose lines are $(k + 1)$-element subsets of $X$, incidence being inclusion. Obviously,
$G_k(N)$ is a $\left({N \choose k}_{N-k}, {N \choose
k+1}_{k+1}\right)$-configuration; hence, we have another $(28_6, 56_3)$-configuration, $G_2(8)$. 

It is straightforward to see that the two $(28_6, 56_3)$-configurations are, in fact, isomorphic. To this end, one simply employs the following bijection between the 28 off-quadric points and the 28 points of $G_2(8)$ (here, by a slight abuse of notation, $X = \{1,2,3,4,5,6,7,8\}$)

\begin{center}
\begin{tabular}{||c|c|c|c||}
\hline \hline
 off-$\mathcal{Q}$ &  $G_2(8)$ & off-$\mathcal{Q}$ &  $G_2(8)$ \\
\hline
~1~  &  $\{1,4\}$  &  ~15~ &  $\{2,3\}$  \\
~2~  &  $\{3,5\}$  &  ~16~ &  $\{4,7\}$ \\
~3~  &  $\{2,5\}$  &  ~17~ &  $\{5,6\}$    \\
~4~  &  $\{4,6\}$  &  ~18~ &  $\{1,5\}$   \\
~5~  &  $\{2,6\}$  &  ~19~ &  $\{1,7\}$    \\
~6~  &  $\{3,6\}$  &  ~20~ &  $\{6,7\}$   \\
~7~  &  $\{1,2\}$  &  ~21~ &  $\{4,5\}$  \\
~8~  &  $\{1,3\}$  &  ~22~ &  $\{7,8\}$  \\
~9~  &  $\{1,6\}$  &  ~23~ &  $\{5,8\}$     \\
~10~ &  $\{2,4\}$  &  ~24~ &  $\{3,8\}$   \\
~11~ &  $\{3,4\}$  &  ~25~ &  $\{2,8\}$    \\
~12~ &  $\{5,7\}$  &  ~26~ &  $\{4,8\}$   \\
~13~ &  $\{3,7\}$  &  ~27~ &  $\{1,8\}$     \\
~14~ &  $\{2,7\}$  &  ~28~ &  $\{6,8\}$     \\
 \hline \hline
\end{tabular}
\end{center}
and verifies step by step that each of the above-listed 56 lines of PG(5,\,2) is also a line of  $G_2(8)$; thus, 
line 1 of PG(5,\,2) corresponds to the line $\{1,4,6\}$ of $G_2(8)$, line 2 to the line $\{1,2,4\}$, line 3 to $\{1,3,4\}$, etc.

This isomorphism entails a very interesting property related to so-called Conwell heptads \cite{con}. Given a $\mathcal{Q}^{+}(5,2)$ of PG(5,\,2), a Conwell heptad is
a set of seven off-quadric points such that each line joining two distinct points of the heptad is skew to the $\mathcal{Q}^{+}(5,2)$. There are altogether eight such heptads: any two of them have a unique point
in common and each of the 28 points off the quadric is contained in two heptads.
The points in the first table are arranged in such a way that, as it is obvious from the
bottom part of the second table, the last seven of them represent a Conwell heptad. From the last table we read off that this particular heptad corresponds to those seven points of $G_2(8)$ whose representatives have mark  `8' in common. Clearly, the remaining seven heptads correspond to those septuples of points of $G_2(8)$ that share one of the remaining seven marks each. Finally, we observe that by removing from our off-quadric $(28_6, 56_3)$-configuration the seven points of a Conwell heptad and all the 21 lines defined by pairs of them one gets a $(21_5, 35_3)$-configuration isomorphic to $G_2(7)$; gradual removal of additional heptads and the corresponding lines yields the following nested sequence of configurations:

\begin{center}
\begin{tabular}{||c|c|c|c||}
\hline \hline
 $\#$ of heptads &  Configuration   & CG       & Remark\\
  removed        &                  &          &         \\
\hline
0                &  $(28_6, 56_3)$  & $G_2(8)$ & \\
1                &  $(21_5, 35_3)$  & $G_2(7)$ & \\
2                &  $(15_4, 20_3)$  & $G_2(6)$ & Cayley-Salmon \\
3                &  $(10_3, 10_3)$  & $G_2(5)$ & Desargues \\
4                &  $(6_2, 4_3)$    & $G_2(4)$ & Pasch\\
5                &  $(3_1, 1_3)$    & $G_2(3)$ & single line \\
6                &  $(1_0, 0_3)$    & $G_2(2)$ & single point\\
7                &                  &          & empty set\\
 \hline \hline
\end{tabular}
\end{center}
Remarkably, this nested sequence of binomial configurations is identical with part of that found to be associated with Cayley-Dickson algebras \cite{cd}.

\section*{Acknowledgment}
This work was partially supported by the VEGA Grant Agency, Project 2/0003/13, as well as by the Austrian Science Fund (Fonds zur F\"orderung der Wissenschaftlichen Forschung (FWF)), Research Project M1564--N27.

\end{document}